\documentclass[11pt]{amsart}

\usepackage[T1]{fontenc}
\usepackage[utf8]{inputenc}
\usepackage{lmodern}
\usepackage{amsmath,amssymb,amsthm}
\usepackage[a4paper,margin=1in]{geometry}
\usepackage{hyperref}

\hypersetup{
  colorlinks=true,
  linkcolor=blue,
  citecolor=blue,
  urlcolor=blue
}

\newtheorem{theorem}{Theorem}
\newtheorem{corollary}{Corollary}

\title[Unbounded Coefficients of Powers]{Unboundedness of the Coefficients of Higher Powers of a Unimodular Power Series}
\author{Zhao Shen}
\address{School of Mathematics and Statistics, Central South University, Changsha, Hunan, China}
\email{sz1021@csu.edu.cn}
\date{}
\subjclass[2020]{30B30, 42A16}
\keywords{unimodular power series, coefficient growth, Parseval's identity, Jensen's inequality}

\begin{document}

\begin{abstract}
Let $R(z)=\sum_{n=0}^{\infty} r_n z^n$ be a power series with $|r_n|=1$ for every $n\ge 0$.
We show that for each integer $m\ge 2$, the coefficient sequence of $R(z)^m$ is unbounded.
The proof combines Parseval's identity with Jensen's inequality.
As a consequence, Conjecture~3.9 of Gawron, Miska, and Ulas \cite{gmu} is confirmed.
\end{abstract}

\maketitle

\section{Introduction}

One classical object in combinatorics and number theory is the Prouhet--Thue--Morse generating function
\[
T(x)=\prod_{n=0}^{\infty}(1-x^{2^n}) = \sum_{n=0}^{\infty} a_n x^n,
\]
Its coefficients are given by $a_n=(-1)^{s_2(n)}\in\{-1,1\}$, where $s_2(n)$ denotes the sum of the binary digits of $n$.
Gawron, Miska, and Ulas \cite{gmu} studied the arithmetic properties of the coefficients of $T(x)^m$
and proposed Conjecture~3.9, which asserts that for every integer $m\ge 2$,
the coefficient sequence of $T(x)^m$ is unbounded.
They proved this for $m=3$ and for $m=2^k$, while the general case remained open.

Since the coefficients of $T(x)$ satisfy $|a_n|=1$, it is natural to place the problem in the broader
setting of power series with \emph{unimodular} coefficients.
In this note we prove that if $R(z)=\sum_{n=0}^{\infty} r_n z^n$ satisfies $|r_n|=1$ for all $n\ge 0$,
then the coefficients of $R(z)^m$ are unbounded for every $m\ge 2$.
This yields Conjecture~3.9 as an immediate corollary.

\begin{theorem}
\label{thm:main}
Let $(r_n)_{n=0}^{\infty}$ be a sequence of complex numbers with $|r_n|=1$ for every $n\ge 0$.
Define $R(z)=\sum_{n=0}^{\infty} r_n z^n$, and for a fixed integer $m\ge 2$ write
\[
R(z)^m=\sum_{n=0}^{\infty} f_m(n) z^n.
\]
Then the sequence $(f_m(n))_{n=0}^{\infty}$ is unbounded.
\end{theorem}

\begin{corollary}
For every integer $m\ge 2$, the coefficient sequence of $T(x)^m$ is unbounded.
In particular, Conjecture~3.9 of \cite{gmu} holds.
Moreover, combining Theorem~\ref{thm:main} with Lemma~3.10 of \cite{gmu} shows that
the coefficients of $T(x)^m$ are unbounded both above and below.
\end{corollary}

\begin{proof}[Proof of Theorem~\ref{thm:main}]
Assume for contradiction that there exists $M>0$ such that
\begin{equation}
\label{eq:bounded}
|f_m(n)|\le M \qquad\text{for all } n\ge 0.
\end{equation}

Since $|r_n|=1$, the Cauchy--Hadamard formula gives radius of convergence $1$,
so $R(z)$ and $R(z)^m$ are holomorphic on the open unit disk $\mathbb{D}$.

Fix $r\in(0,1)$ and write $z=re^{i\theta}$. Applying Parseval's identity to
$R(z) = \sum r_n z^n$ on $|z|=r$:
\begin{equation}
\label{eq:parseval1}
\frac{1}{2\pi}\int_0^{2\pi}|R(re^{i\theta})|^2\,d\theta
=\sum_{n=0}^{\infty}|r_n|^2 r^{2n}
=\frac{1}{1-r^2}.
\end{equation}

Applying Parseval's identity to the holomorphic function $R(z)^m = \sum f_m(n) z^n$
on $|z|=r$, and using \eqref{eq:bounded}:
\begin{equation}
\label{eq:parseval2}
\frac{1}{2\pi}\int_0^{2\pi}|R(re^{i\theta})|^{2m}\,d\theta
=\sum_{n=0}^{\infty}|f_m(n)|^2 r^{2n}
\le \frac{M^2}{1-r^2}.
\end{equation}

Since $\phi(x)=x^m$ is convex on $[0,\infty)$ for $m\ge 2$, Jensen's inequality gives
\begin{equation}
\label{eq:lowerbound}
\frac{1}{2\pi}\int_0^{2\pi}|R(re^{i\theta})|^{2m}\,d\theta
\ge \left(\frac{1}{2\pi}\int_0^{2\pi}|R(re^{i\theta})|^2\,d\theta\right)^m
= \left(\frac{1}{1-r^2}\right)^m.
\end{equation}

Combining \eqref{eq:parseval2} and \eqref{eq:lowerbound}:
\[
\frac{1}{(1-r^2)^{m-1}}\le M^2 \qquad (0<r<1).
\]
Since $m\ge 2$, the left-hand side tends to $+\infty$ as $r\to 1^{-}$, a contradiction.
\end{proof}

\begin{proof}[Proof of Corollary]
The coefficients of $T(x)$ lie in $\{-1,1\}$, so Theorem~\ref{thm:main} applies directly.
The final statement follows from Lemma~3.10 of \cite{gmu}.
\end{proof}

\section*{Acknowledgements}
The author thanks Luo Sijie, Liu Sihan, and Liu Zhecheng for helpful comments and suggestions.
This work was supported by the National Natural Science Foundation of China (No.~12301011).

\end{document}